\documentclass[a4paper,11pt]{article}
\usepackage{graphicx}
\usepackage{a4wide}
\usepackage[english]{babel}
\usepackage{amsfonts}
\usepackage{amsmath}
\usepackage{amssymb}
\usepackage{dsfont}
\usepackage{amsthm}
\newtheorem{lemma}{Lemma}
\newtheorem{theorem}{Theorem}

\usepackage{float}

\newcommand {\E} {\mathbb{E}}

\newcommand {\p} {\mathbb{P}}
\newcommand {\Z} {\mathbb{Z}}
\newcommand {\N} {\mathbb{N}}

\newcommand {\R} {\mathbb{R}}

\newcommand {\ve} {\varepsilon}

\newcommand {\no} {\mathcal{N}}
\newcommand {\nd} {\mathcal{N}^{(d)}}

\makeatletter\def\blfootnote{\xdef\@thefnmark{}\@footnotetext}\makeatother
\allowdisplaybreaks
\title{\bf On the limit distribution of the normality measure of random binary sequences}
\author{Christoph Aistleitner\footnote{Department of Applied Mathematics, School of Mathematics and Statistics, University of New South Wales, Sydney NSW 2052, Australia. \mbox{e-mail}: \texttt{aistleitner@math.tugraz.at}. Research supported by a Schr\"odinger scholarship of the Austrian Research
Foundation (FWF).}}
\begin{document}

\date{}
\maketitle

\blfootnote{{\bf MSC 2010:} Primary Classification: 68R15, Secondary Classification: 11K45, 60C05, 60F05}
\blfootnote{{\bf keywords:} random sequence, pseudorandom sequence, normality measure, discrepancy, limit distribution}

\begin{abstract}
We prove the existence of a limit distribution for the normalized normality measure $\mathcal{N}(E_N)/\sqrt{N}$ (as $N \to \infty$) for random binary sequences $E_N$, by this means confirming a conjecture of Alon, Kohayakawa, Mauduit, Moreira and R{\"o}dl. The key point of the proof is to approximate the distribution of the normality measure by the exiting probabilities of a multidimensional Wiener process from a certain polytope.
\end{abstract}

\section{Introduction and statement of results}

In a series of papers starting in 1997, Mauduit and S{\'a}rk{\"o}zy~\cite{ms1} introduced and studied several measures of pseudorandomness for finite binary sequences. In the present paper we will mainly be concerned with the \emph{normality measure} $\mathcal{N}(E_N)$; however, for comparison the connection between our new results in the present paper and earlier results for other pseudorandomness measures will be described in Section~\ref{sec2} below. Let a finite binary sequence $E_N = (e_1, \dots, e_N) \in \{-1,1\}^N$ be given. For $k \in \N, ~M \in \N$ and $X \in \{-1,1\}^k$, set 
$$
T(E_N,M,X) = \# \left\{ n:~0 \leq n < M, ~\textrm{and}~(e_{n+1},\dots,e_{n+k}) = X \right\}.
$$
Thus $T(E_N,M,X)$ counts the number of occurrences of the pattern $X$ among the first $M+k$ elements of $E_N$. The normality measure $\mathcal{N}(E_N)$ is defined as
\begin{equation} \label{defno}
\mathcal{N} (E_N) = \max_k ~\max_X ~\max_M \left| T(E_N,M,X) - \frac{M}{2^k}\right|,
\end{equation}
where the maxima are taken over $k \leq \log_2 N,~X \in \{-1,1\}^k,~1 \leq M \leq N+1-k$.\\ 

Alon, Kohayakawa, Mauduit, Moreira and R{\"o}dl studied the \emph{minimal}~\cite{akmmr_min} and \emph{typical}~\cite{akmmr_typ} values of $\mathcal{N}(E_N)$. Concerning the minimal possible value of $\mathcal{N}(E_N)$, they proved
\begin{equation} \label{gap}
\left( \frac{1}{2} - \ve \right) \log_2 N \leq \min_{E_N \in \{-1,1\}^N} \mathcal{N}(E_N) \leq 3 N^{1/3} (\log N)^{2/3}
\end{equation}
for sufficiently large $N$. There is a relatively large gap between the lower and upper bound in \eqref{gap}, but apparently no further progress has been made since the publication of~\cite{akmmr_min} in 2006. Thus the question asking for the minimal possible order of $\mathcal{N}(E_N)$ remains unsolved.\footnote{Since the initial submission of the present manuscript, this problem has been solved up to logarithmic factors. See \cite{ai2} for details.} For constructions of sequences having small normality measure, see also~\cite{rivat}.\\

Concerning the typical value of $\mathcal{N}(E_N)$, Alon \emph{et al.} showed (improving earlier results of Cassaigne, Mauduit and S{\'a}rk{\"o}zy~\cite{cms}) that choosing $E_N$ randomly from $\{-1,1\}^N$, for any $\ve>0$ there exist $\delta_1,\delta_2>0$ (depending on $\ve$) such that 
\begin{equation} \label{n1}
\delta_1 \sqrt{N} < \mathcal{N} (E_N) < \delta_2 \sqrt{N}
\end{equation}
with probability at least $1 - \ve$ for sufficiently large $N$.  The lower bound is optimal, since for any $\delta>0$
\begin{equation} \label{n2}
\liminf_{N \to \infty} ~\p \left(\no(E_N) < \delta \sqrt{N}\right) > 0.
\end{equation}
Summarized, the normality measure of a random binary sequence is typically of order $\sqrt{N}$. At the end of~\cite{akmmr_typ}, Alon \emph{et al.} pose the problem of investigating the existence of a limit distribution of 
\begin{equation} \label{lex}
\frac{\mathcal{N}(E_N)}{\sqrt{N}}
\end{equation}
for random $E_N$ as $N \to \infty$, and write that it is ``most likely'' that such a limit distribution exists. The purpose of the present paper is to confirm their conjecture, and show that a limit distribution of \eqref{lex} in fact exists.

\begin{theorem} \label{th1}
There exists a probability distribution function $F$ such that for random $E_N$ the distribution of $\mathcal{N}(E_N)/\sqrt{N}$ converges to $F$ as $N \to \infty$. More precisely, for $E_N$ having uniform distribution in $\{-1,1\}^N$ for $N \geq 1$, we have for any $t \in \R$ that
$$
\lim_{N \to \infty} \p \left( \frac{\mathcal{N}(E_N)}{\sqrt{N}} \leq t \right) = F(t).
$$
The function $F(t)$ is continuous for all $t \in \R$.
\end{theorem}

\section{Comparison with other pseudorandomness measures} \label{sec2}

Two other measures of pseudorandomness, introduced by Mauduit and S{\'a}rk{\"o}zy, are the \emph{well-distribution measure} $W(E_N)$ and the \emph{correlation measure} $C_k(E_N)$. For $M \in \N,~a \in \Z$ and $b \in \N$ set
$$
U (E_N,M,a,b) = \sum \big\{ e_{a+jb}:~1 \leq j \leq M,~1 \leq a + jb \leq N~\textrm{for all} ~j \big\}.
$$
The well-distribution measure $W(E_N)$ is then defined as
$$
W(E_N) := \max \left\{ \left| U(E_N,M,a,b)\right|, ~\textrm{where}~1 \leq a+b \textrm{~and~}  a+Mb \leq N \right\},
$$
hence measuring the maximal discrepancy of $E_N$ along an arithmetic progression. This well-distribution measure can be seen as a special case of a combinatorial discrepancy measure (cf.~\cite[Chapter 4]{mat}), but can also be directly modified into a generalization of the concept of discrepancy in the context of uniform distribution modulo one in analytic number theory~(see~\cite{bpt}).\\

For $k \in \N, ~M \in \N$ and $D=(d_1, \dots, d_k) \in \N^k$ with $0 \leq d_1 < \dots < d_k < N$ we define
\begin{eqnarray*}
V(E_N,M,D) & = & \sum \left\{ e_{n+d_1} \dots e_{n+d_k}:~1 \leq n \leq M,~n + d_k \leq N \right\}.
\end{eqnarray*}
Thus $V(E_N,M,D)$ measures the correlation among $k$ segments of $E_N$, which are relatively positioned according to $D$. The correlation measure of order $k$, which is denoted by $C_k(E_N)$, is defined as 
$$
C_k (E_N) = \max \left\{ | V(E_N,M,D)|:~M,D \textrm{~satisfy~} M + d_k \leq N \right\}.
$$ 
Note that contrary to the normality measure $\mathcal{N}(E_N)$ and the well-distribution measure $W(E_N)$, which only depend on the sequence $E_N$, the correlation measure $C_k(E_N)$ depends on an additional parameter $k$.\\

In~\cite{akmmr_typ} several results concerning the typical asymptotic order of the well-distribution measure and the correlation measure are proved. For random $E_N$, for any $\ve>0$ and sufficiently large $N$ the correlation measure satisfies
$$
\frac{2}{5} \sqrt{N \log \binom{N}{k}} \leq C_k(E_N) \leq \frac{7}{4} \sqrt{N \log \binom{N}{k}}
$$
with probability at least $1-\ve$, for all $2 \leq k \leq N/4$. Furthermore, for any function $k(N) \leq \log N - \log \log N$ with probability at least $1 - \ve$ we have
$$
1 - \ve < \frac{C_k(E_N)}{\E(C_k)} < 1 + \ve,
$$
for sufficiently large $N$.\footnote{Note added in proof: a strong improvement of these results has recently been announced by Schmidt. See~\cite{schmidt}.} In other words, the correlation measure is concentrated around its mean if $k$ is ``small''; consequently, the limit distribution is in this case the Dirac measure centered at $\E(C_k)$. For recent results on the correlation measure, particularly concerning its dependence on the parameter $k$, see~\cite{an,g1,g2}.\\

For the well-distribution measure, in~\cite{akmmr_typ} the following results are proved: for random $E_N$, for any $\ve>0$ there exist $\delta_1, \delta_2 > 0$ such that
$$
\delta_1 \sqrt{N} < W(E_N) < \delta_2 \sqrt{N}
$$
with probability at least $1-\ve$, for $N$ sufficiently large. The lower bound is optimal, since for any $\delta>0$
$$
\liminf_{N \to \infty} ~\p \left(W(E_N) < \delta \sqrt{N}\right) > 0;
$$
note that these two results for the well-distribution measure are similar to those for the normality measure in \eqref{n1} and \eqref{n2}. For other recent results concerning the well-distribution measure, see e.g.~\cite{mer, sz, wu}. The existence of a limit distribution of the normalized well-distribution measure $W_N(E_N)/\sqrt{N}$ of a random sequence was conjectured in~\cite{akmmr_typ}, and confirmed in~\cite{ai}. In~\cite{ai} I wrote 
\begin{quote} 
``The case of the normality measure $\no(E_N)$ seems to be much more difficult''
\end{quote} 
and this is really the case. The problem is mainly caused by the fact that patterns can overlap, and the independence of the occurrences of the single digits (a \emph{digit} meaning here one of the numbers $-1$ and $1$) is not sufficient to deduce independence of the occurrences of (possibly overlapping) patterns of multiple digits; for example, the pattern \emph{11} appears in the block of digits \emph{111} twice, and if in a block of three random digits $e_1 e_2 e_3$ the first two digits are \emph{11} there is an increased probability that the pattern \emph{11} will also appear in the last two digits (namely, the probability is 1/2, while it should be 1/4 in the independent case).\\
We overcome the problem by cutting the index set $\{1, \dots, N\}$ into blocks of length $d$ (for large $d$) and initially considering only the appearances of patterns entirely contained \emph{within} one of these blocks of digits. The occurrence of a certain block $B$ (out of $2^d$ possible block) can be interpreted as the action of a random walk on a $2^d$-dimensional lattice (moving one step into the direction associated with this specific block); the probabilities of the normality measure exceeding a certain value are then asymptotically equal to the probabilities of a corresponding Wiener process (the limit process of the normalized random walk) leaving a certain polytope. Furthermore, we use decorrelation methods to show that the impact of the occurrences of patterns not entirely contained within a block of length $d$ is small (for $d$ sufficiently large).

\section{Auxiliary results}


\begin{lemma}[{Maximal Bernstein inequality; see e.g.~\cite[Lemma
2.2]{einmahl}}] \label{bernstein}
For a sequence $\xi_1, \dots, \xi_N$ of independent, identically distributed (i.i.d.) random variables having mean zero and
variance $\sigma^2$, and satisfying $|\xi_n| \leq 1$, we have for $t \geq 0$
$$
\p \left( \max_{1 \leq M \leq N} \left| \sum_{n=1}^M \xi_n \right| > t  \right)
\leq 2 e^{-t^2/(2 N \sigma^2 + 2t/3)}
$$
\end{lemma}

The following lemma is a special case of the multidimensional version of Donsker's theorem; see e.g. {\cite[Theorem 4.3.5]{ww}}.
\begin{lemma} \label{lemmadon}
Let $(\xi_n)_{n \geq 1}$ be a sequence of bounded, i.i.d., $d$-dimensional random vectors having expectation zero and covariance matrix $\Sigma$. For $s \in [0,1]$, set
$$
Y_N(s) = \frac{1}{\sqrt{N}} \sum_{n=1}^{\lceil Ns \rceil} \xi_n.
$$
Then 
$$
Y_N \Rightarrow Z,
$$
where $Z$ is a $d$-dimensional Wiener process with mean zero and covariance matrix $\Sigma$, and $\Rightarrow$ denotes weak convergence in the $d$-dimensional Skorokhod space $D([0,1]^d)$. $Z$ satisfies the equation $Z = B C$, where $B$ is a $d$-dimensional standard Wiener process (standard Brownian motion)  and $C$ is a $d \times d$ matrix for which $C^T C = \Sigma$ holds.
\end{lemma}

\emph{Remark:~} The use of deep probabilistic results such as Donsker's theorem could be avoided by replacing the continuous parameter $s \in [0,1]$ by taking $s$ from a discrete, equidistant grid $\{s_0, \dots, s_R\}$, which corresponds to approximating the normality measure by a variant for which only the values $s_0 N, \dots, s_R N$ are allowed for $M$. Then we could prove Theorem~\ref{th1} using the classical multidimensional central limit theorem for i.i.d.\ random vectors, and choosing $R$ ``large''. However, the proof is much shorter and clearer using Donsker's theorem.

\section{Idea of the proof and preliminaries}

For the rest of this paper, we assume that $(e_n)_{n \geq 1}$ are i.i.d.\ random variables, each taking a value from $\{-1,1\}$ equiprobably. We write $E_N = (e_1, \dots, e_N)$. We will split the proof of Theorem~\ref{th1} into several parts. The main idea of the proof is to show that the value of the normality measure $\mathcal{N}$ is with large probability almost equal to a restricted normality measure $\mathcal{N}^{(d)}$ considering only those occurrences of an pattern $X$ for which the index set is entirely contained in a set of the form $\{md+1, \dots, (m+1)d\}$ for some $m$. For each $m$ the block $(e_{md+1}, \dots, e_{(m+1)d})$ is a random element of $\{-1,1\}^d$, and of course for different values of $m$ these blocks are independent.\\
To each possible block $B \in \{-1,1\}^d$ of $d$ digits we will assign a vector of the standard (Cartesian) base of $\R^{2^d}$, and interpret the occurrence of a specific block $B$ as a step forward of a random walk in the direction of the assigned coordinate (subtracting the expected value, that is the average of all possible movements, in each step). By Donsker's theorem this random walk (appropriately normalized) converges weakly to a Wiener process.\\

The probability of this restricted normality measure exceeding a certain value then converges to the exiting probability of the limiting Wiener process from a certain polytope. This may sound surprising, but in fact is quite natural if one reflects on when the normality measure exceeds a certain value: the restricted normality measure $\mathcal{N}^{(d)} (E_N)$ is larger than $t$ if for some possible pattern $X$ the sum (over $B \in \{-1,1\}^d$) of the number of occurrences of $B$ multiplied with the number of occurrences of $X$ within $B$ is greater than $t$.\\
This sum of products can be written as a scalar product of the position of the random walk, multiplied with a ``weight'' vector (counting the number of occurrences of $X$ within $B$ for each possible $B$), and the event of this scalar product exceeding a certain value $t$ equals the event of the random walk exiting the space between two hyperplanes (depending on $t$). Taking the maximum over all possible values of $X$ in the definition of the normality measure corresponds to taking the exiting probabilities of the random walk from the intersection of all possible hyperplanes associated to some $X$ (which in our case in fact produces a proper polytope).\\

Let $d \geq 1$ be given, and for $m \geq 0$ set
$$
\Delta_m = \{md+1, \dots, (m+1)d\}.
$$
Modifying the definition of $\no(E_N)$ in \eqref{defno}, we define a restricted normality measure $\nd(E_N)$ by setting
\begin{eqnarray*}
& & T^{(d)}(E_N,M,X) \\
& = & \# \left\{ n:~0 \leq n < M, ~(n ~\textup{mod}~ d) \in \{0, \dots, d-k\}, ~\textrm{and}~(e_{n+1},\dots,e_{n+k}) = X \right\}
\end{eqnarray*}
and
\begin{equation} \label{ris}
\nd (E_N) = \max_k ~\max_X ~\max_M \left|  T^{(d)}(E_N,M,X) - \frac{M+k-1}{d} \left(\frac{d-k+1}{2^k} \right)\right|,
\end{equation}
where the maxima are taken over $k \leq d,~X \in \{-1,1\}^k,~1 \leq M \leq N+1-k$ and $M+k-1 \equiv 0 \mod d$. That means that for $\nd (E_N)$ we consider only patterns $X$ of length at most $d$, and only those occurrences $(e_{n+1},\dots,e_{n+k})$ of such patterns for which the index set $\{n+1,\dots,n+k\}$ is entirely contained in a set $\Delta_m$ for some $m \geq 0$ (that means, the index set $\{n+1,\dots,n+k\}$ does not overstep any integer multiple of $d$). The additional assumption $M+k-1 \equiv 0 \mod d$ means that we only consider sets of indices $\{1, \dots, M+k\}$ which end at an integer multiple of $d$, and can therefore be written as a union of sets $\Delta_m$; considering only values $M$ of this form accounts for an additional possible error of at most $d$ in comparison with a normality measure defined without this additional restriction.\\
Note that the number of possible index sets $\{n+1, \dots, n+k\}$ satisfying $n < M$ for $M+k-1 \equiv 0 \mod d$, which do \emph{not} overstep any integer multiple of $d$, is precisely $(M+k-1)(d-k+1)/d$ (which accounts for the very last term in line \eqref{ris}). In contrast, the number of such index sets which \emph{do} overstep an integer multiple of $d$ is $\left(\frac{M+k-1}{d} -1 \right)(k-1)$ (which accounts for the very last term in line \eqref{ter1} below). Naturally,  
$$
\frac{(M+k-1) (d-k+1)}{d} + \left(\frac{M+k-1}{d} -1 \right)(k-1) = M.
$$

Furthermore, we set
\begin{eqnarray*}
& & \bar{T}^{(d)}(E_N,M,X) \\
& = & \# \left\{ n:~0 \leq n < M, ~(n ~\textup{mod}~ d) \in \{d-k+1, \dots, d-1\}, ~\textrm{and}~(e_{n+1},\dots,e_{n+k}) = X \right\}.
\end{eqnarray*}
Thus $\bar{T}^{(d)}(E_N,M,X)$ counts the number of occurrences of a pattern $X$ for which the index set is \emph{not} entirely contained in $\Delta_m$ for some $m$; clearly this means that 
$$
T(E_N,M,X) = T^{(d)}(E_N,M,X) + \bar{T}^{(d)}(E_N,M,X).
$$
Consequently, $\no (E_N)$ is bounded above by
\begin{align}
& \max\left\{ \nd (E_N) + d + \max_{k \leq d} ~\max_X ~\max_M \left| \bar{T}^{(d)}(E_N,M,X) - \left(\frac{M+k-1}{d} -1 \right)\frac{k-1}{2^k}\right|,\right. \label{ter1}\\
 & \left. \qquad \quad \max_{d < k \leq \log_2 N} ~\max_X ~\max_M \left| T(E_N,M,X) - \frac{M}{2^k} \right| \right\},\label{ter2}
\end{align}
where in each line the maxima are taken over $X \in \{-1,1\}^k$ and $1 \leq M \leq N+1-k$, respectively, and in the last maximum in line \eqref{ter1} we have the additional restriction $M+k-1 \equiv 0 \mod d$. Here, to compare $\nd$ with $\no$, the additional term $d$ in line \eqref{ter1} comes from the restriction $M+k-1 \equiv 0 \mod d$ in the definition of $\nd(E_N)$, and the last term in line \eqref{ter1} adds the contribution of the occurrences of a pattern for which the index set is not entirely contained in $\Delta_m$ for some $m$. The term in line \eqref{ter2} accounts for the possible contribution of patterns which consist of more than $d$ digits. As we will see the contribution of the last term in \eqref{ter1} and of the term in \eqref{ter2} is with large probability very small, provided $d$ is chosen sufficiently large.\\

On the other hand, $\no (E_N)$ is bounded below by
\begin{align}
 & \nd (E_N) - \max_{k \leq d} ~\max_X ~\max_M \left| \bar{T}^{(d)} (E_N,M,X) - \left(\frac{M+k-1}{d} -1 \right)\frac{k-1}{2^k} \right|, \label{ter3}
\end{align}
where the second and the third maximum are taken over $X \in \{-1,1\}^k$ and $1 \leq M \leq N+1-k$, respectively.\\

Let $\mathcal{B}$ denote the set $\{-1,1\}^d$ of all blocks of $d$ digits, and let $\mathcal{X}$ denote the set of all possible patterns of at most $d$ digits, that is all patterns $X \in \bigcup _{k \leq d} \{-1,1\}^k$. Furthermore, for any $B \in \mathcal{B}$ and $X \in \mathcal{X}$, let $w_{X,B}$ denote the number of occurrences of the pattern $X$ within the block $B$, that is
$$
w_{X,B} = T(B,d-k+1,X).
$$
In the sequel, we will use the symbol ``$\cdot$'' for the scalar product of two vectors, and for the product of a scalar and a vector. We will generally write vectors in bold font.\\

We enumerate the $2^d$ elements of $\mathcal{B}$ by $b_1, \dots, b_{2^d}$, and let $\boldsymbol{\beta}_1, \dots, \boldsymbol{\beta}_{2^d}$ denote the Cartesian base vectors of $\R^{2^d}$. Furthermore, writing $\mathbf{1}$ for the $d$-dimensional vector $(1, \dots, 1)$, we set
\begin{equation} \label{xm}
X_m = \left( \sum_{u=1}^{2^d} \mathds{1}_{b_u}(e_{md+1},\dots,e_{(m+1)d}) \cdot \boldsymbol{\beta}_u \right) - \frac{1}{2^d} \cdot \mathbf{1} .
\end{equation}
Here $\mathds{1}_{b_u}(e_{md+1},\dots,e_{(m+1)d})$ is the indicator function of $b_u$, meaning $\mathds{1}_{b_u}(e_{md+1},\dots,e_{(m+1)d})=1$ if $(e_{md+1},\dots,e_{(m+1)d})=b_u$, and $\mathds{1}_{b_u}(e_{md+1},\dots,e_{(m+1)d})=0$ otherwise. Then $X_m$ is a random vector for $m \geq 0$, and $X_m = \boldsymbol{\beta}_u - 2^{-d} \cdot \mathbf{1}$ if and only if $(e_{md+1},\dots,e_{(m+1)d}) = b_u$. The random vectors $X_m, ~m \geq 1,$ are independent, identically distributed, and have expectation zero.\\

For any $X \in \mathcal{X}$ we set $\mathbf{w}_X = (w_{X,b_1}, \dots, w_{X,b_u}) \in \R^{2^d}$, 
and for $t \geq 0$ define a polytope $P^{(t)}$ as
$$
P^{(t)} = \bigcap_{X \in \mathcal{X}} \left\{ \mathbf{y} \in \R^{2^d}:~\left|\mathbf{w}_X \cdot \mathbf{y} \right| \leq t \right\}.
$$
Note that $P^{(t)}$ is defined as the intersection of finitely many half-spaces; however, it is in fact a proper polytope, since for all $\ell$ we have $b_\ell \in \mathcal{X}$, and, since $w_{b_\ell,b_u} = 1$ if and only if $\ell = u$,
$$
\mathbf{w}_{b_\ell} = \sum_{u=1}^{2^d} w_{b_\ell,b_u} \boldsymbol{\boldsymbol{\beta}}_u = \boldsymbol{\boldsymbol{\beta}}_\ell.
$$
Consequently the absolute value of the $\ell$-th coordinate of any element $\mathbf{y} \in P^{(t)}$ is bounded by $t$, for $\ell \in \{1, \dots, 2^d\}$. In other words, $P^{(t)}$ is contained in the $2^d$-dimensional cube $[-t,t]^{2^d}$. Furthermore, $P^{(t)}$ is nonempty if $t > 0$.\\

\section{Main lemmas}

Lemma~\ref{lemmacru} below is the key ingredient in the proof, showing that the probability of the normalized restricted normality measure $\mathcal{N}^{(d)}(E_N)/\sqrt{N}$ exceeding the value $t$ converges to the exiting probabilities of an appropriate Wiener process from the polytope $P^{(t)}$. The following Lemmas~\ref{lemma1} and~\ref{lemma2} state that the error made by approximating the normality measure $\no$ by the restricted normality measure $\nd$ is ``small'', provided $d$ is ``large''. 

\begin{lemma} \label{lemmacru}
Let $\Sigma$ be the covariance matrix of the random variables $X_m$ in \eqref{xm}. Then for any $t \in \R$,
$$
\lim_{N \to \infty} \p \left( \frac{\nd(E_N)}{\sqrt{N}} \leq t \right) = \p \left( \frac{Z(s)}{\sqrt{d}} \in P^{(t)} ~\textup{for all~} s \in [0,1] \right),
$$
where $Z$ is a $d$-dimensional Wiener process with mean zero and covariance matrix $\Sigma$.
\end{lemma}

\begin{lemma} \label{lemma1}
For any $d \geq	4$ and $N \geq N_0(d)$
\begin{align*}
& \p \left( \max_{k \leq d} ~\max_{X \in \{-1,1\}^k}  ~\max_{\substack{1 \leq M \leq N+1-k,\\M+k-1 \equiv 0 \mod d}} \left| \bar{T}^{(d)}(E_N,M,X) \right.\right. \\
& \left.\left. \qquad \qquad \qquad - \left(\frac{M+k-1}{d} -1 \right)\frac{k-1}{2^k} \right| \geq \frac{6 \sqrt{N} \sqrt{\log d}}{\sqrt{d}} \right) < \frac{1}{d^2-1}.
\end{align*}
\end{lemma}

\begin{lemma} \label{lemma2}
For any $d \geq 4$ and $N \geq N_0(d)$
$$
\p \left( \max_{d < k \leq \log_2 N} ~\max_{X \in \{-1,1\}^k} ~\max_{1 \leq M \leq N + 1 -k} \left| T(E_N,M,X) - \frac{M}{2^k} \right| \geq \frac{16 \sqrt{N}}{d} \right) < \frac{1}{d^{2d}}.
$$
\end{lemma}

\emph{Proof of Lemma~\ref{lemmacru}:~} The main part of the proof of Lemma~\ref{lemmacru} is to identify the probabilities 
$\p \left( \frac{\nd(E_N)}{\sqrt{N}} > t \right)$ as the exiting probabilities of the random walk 
$$
\frac{\sum_{m=1}^{\lceil Ns/d \rceil} X_m}{\sqrt{N}}, \qquad s \in [0,1],
$$
from the polytope $P^{(t)}$. The convergence is then an immediate consequence of Donsker's theorem in Lemma~\ref{lemmadon}.\\

By the restriction on the values of $M$ in the definition of $\nd(E_N)$, we can assume that $N$ is an integer multiple of $d$, which means that $N = Rd$ for some $R$. Now
\begin{equation} \label{equa}
\nd (E_N) > t \sqrt{N}
\end{equation}
if and only if there exists an pattern $X \in \mathcal{X}$ of $k \leq d$ digits and an $M \leq N+1-k$ satisfying $M+k-1 \equiv 0 \mod d$ such that
$$
\left|  T^{(d)}(E_N,M,X) - \frac{1}{2^k} \frac{(M+k-1)(d-k+1)}{d} \right| > t \sqrt{N},
$$
which happens if and only if
\begin{equation} \label{equi1}
\left|  \underbrace{\left( \sum_{m=1}^{(M+k-1)/d} \mathds{1}_{b_u} (e_{(m-1)d+1},\dots,e_{md}) ~ w_{X,b_u}  \right)}_{=\sum_{m=1}^{(M+k-1)/d} \left( X_m + 2^{-d} \cdot \mathbf{1} \right) \cdot \mathbf{w}_X} - \frac{1}{2^k} \frac{(M+k-1)(d-k+1)}{d} \right|  > t \sqrt{N}.
\end{equation}
Some simple combinatorics shows that
$$
\sum_{u=1}^{2^d} w_{X,b_u} = \mathbf{w}_X \cdot \mathbf{1} = 2^{d-k} (d-k+1),
$$
and thus \eqref{equi1} is equivalent to
$$
\frac{\left| \left(\sum_{m=1}^{(M+k-1)/d} X_m\right) \cdot \mathbf{w}_X \right|}{\sqrt{N}} > t,
$$
which furthermore is equivalent to
$$
\max_{1 \leq s \leq R} \frac{\left| \left(\sum_{m=1}^{sN/d} X_m\right) \cdot \mathbf{w}_X \right|}{\sqrt{N}} > t.
$$
Thus \eqref{equa} is equivalent to
$$
\frac{\sum_{m=1}^{sN/d} X_m}{\sqrt{N}} \not\in P^{(t)} \qquad \textrm{for some} \qquad s \in \left\{\frac{1}{R},\dots,\frac{R}{R} \right\},
$$
and consequently also equivalent to
$$
\frac{\sum_{m=1}^{\lceil sN/d \rceil} X_m}{\sqrt{N}} \not\in P^{(t)} \qquad \textrm{for some} \qquad s \in [0,1],
$$
By Donsker's theorem this implies that
$$
\lim_{N \to \infty} \p \left(\nd(E_N) > t \sqrt{N} \right) = \p \left( \frac{Z(s)}{\sqrt{d}} \not\in P^{(t)} \quad \textup{for some} \quad s \in [0,1] \right),
$$
which proves the lemma. \qquad $\square$\\

\emph{Proof of Lemma~\ref{lemma1}:~} We can again assume that $N$ is divisible by $d$. Clearly, 
\begin{align}
& \p \left( \max_{k \leq d} ~\max_{X \in \{-1,1\}^k} ~\max_M \left| \bar{T}(E_N,M,X) - \left(\frac{M+k-1}{d} -1 \right)\frac{k-1}{2^k} \right| \geq \frac{6 \sqrt{N} \sqrt{\log d}}{\sqrt{d}} \right) \nonumber\\
\leq & \sum_{k \leq d} ~\sum_{X \in \{-1,1\}^k} \p \left(\max_M  \left| \bar{T}(E_N,M,X) - \left(\frac{M+k-1}{d} -1 \right)\frac{k-1}{2^k} \right| \geq \frac{6 \sqrt{N} \sqrt{\log d}}{\sqrt{d}} \right), \label{fre}
\end{align}
where in both lines the last maximum is taken over $1 \leq M \leq N+1-k,~M+k-1 \equiv 0 \mod d$.\\

Let $k \leq d$ and $X \in \{-1,1\}^k$ be fixed. We introduce the notation
\begin{eqnarray} 
\bar{T}^{(d)}_r (E_N,M,X) & = & \# \left\{ n_:~0 \leq n < M, ~n \equiv r ~\textup{mod}~ d, ~\textrm{and}~(e_{n+1},\dots,e_{n+k}) = X \right\}.\label{not}
\end{eqnarray}
Then 
$$
\bar{T}^{(d)} (E_N,M,X) = \sum_{r=d-k+1}^{d-1} \bar{T}^{(d)}_r (E_N,M,X),
$$
and consequently
\begin{eqnarray*}
& & \p \left(\max_{\substack{1 \leq M \leq N+1-k,\\M+k-1 \equiv 0 \mod d}}  \left| \bar{T}(E_N,M,X) - \frac{1}{2^k} \frac{(M+k-1)(k-1)}{d} \right| \geq \frac{6 \sqrt{N} \sqrt{\log d}}{\sqrt{d}} \right) \\
& \leq & \sum_{r=d-k+1}^{d-1} \p \left(\max_{\substack{1 \leq M \leq N+1-k,\\M+k-1 \equiv 0 \mod d}}  \left| \bar{T}^{(d)}_r (E_N,M,X) - \frac{M+k-1}{2^k d} \right| \geq \frac{6 \sqrt{N} \sqrt{\log d}}{k \sqrt{d}} \right).
\end{eqnarray*}
Note that for $M$ satisfying $M+k-1 \equiv 0 \mod d$ we have
$$
\bar{T}^{(d)}_r (E_N,M,X) = \sum_{m=1}^{(M+k-1)/d} \mathds{1}_X \left( e_{(m-1)d+r+1}, \dots, e_{(m-1)d+r+k} \right).
$$
By the independence of $(e_n)_{n \geq 1}$ the random variables $\mathds{1}_X \left( e_{md+r+1}, \dots, e_{md+r+k} \right)$ are also independent for $1 \leq m \leq N/d$ (they are constructed in such a way that the indices do not overlap for different values of $m$). The random variables $\mathds{1}_X \left( e_{md+r+1}, \dots, e_{md+r+k} \right) - 2^{-k}$ have mean zero and variance $2^{-k}(1-2^{-k}) \leq 2^{-k}$. Using Lemma~\ref{bernstein} for $t = \frac{3 \sqrt{k} 2^{-k/2} \sqrt{N} \sqrt{\log d}}{\sqrt{d}}$, we obtain
\begin{eqnarray*}
& & \p \left(\max_{\substack{1 \leq M \leq N+1-k,\\M+k-1 \equiv 0 \mod d}}  \left| \bar{T}^{(d)}_r (E_N,M,X) - \frac{M+k-1}{2^k d} \right| \geq \frac{3 \sqrt{k} 2^{-k/2} \sqrt{N} \sqrt{\log d}}{\sqrt{d}} \right) \\
& \leq & 2 \exp \left(- \frac{9 k 2^{-k} N d^{-1} \log d}{2Nd^{-1}2^{-k} + 6 \sqrt{k} 2^{-k/2} \sqrt{N} \sqrt{\log d} ~d^{-1/2}/3} \right) \\
& \leq & 2 d^{-3k}
\end{eqnarray*}
for sufficiently large $N$. Note that for $k \geq 1$,
\begin{equation*} \label{and}
3 \sqrt{k} 2^{-k/2} \leq \frac{6}{k},
\end{equation*}
and consequently
$$
\p \left(\max_{\substack{1 \leq M \leq N+1-k,\\M+k-1 \equiv 0 \mod d}}   \left| \bar{T}^{(d)}_r (E_N,M,X) - \frac{M}{2^k d} \right| \geq \frac{6 \sqrt{N} \sqrt{\log d}}{k \sqrt{d}} \right) \leq 2 d^{-3k}.
$$

Thus, using \eqref{fre} and the assumption that $d \geq 4$, this implies
\begin{align*}
& \p \left( \max_{k \leq d} ~\max_{X \in \{-1,1\}^k} ~\max_{\substack{1 \leq M \leq N+1-k,\\M+k-1 \equiv 0 \mod d}} \left| \bar{T}(E_N,M,X) - \left(\frac{M+k-1}{d} -1 \right)\frac{k-1}{2^k} \right| \geq \frac{6 \sqrt{N} \sqrt{\log d}}{\sqrt{d}} \right) \nonumber\\
 \leq & \sum_{k \leq d} ~\sum_{X \in \{-1,1\}^k} ~\sum_{r=d-k}^{d-1} 2 d^{-3k} \\
 \leq & 2 \sum_{k \leq d} 2^{k} k d^{-3k} \\
 \leq & \sum_{k \leq d} d^{-2k} \\
 \leq & \frac{1}{d^2-1}. \qquad \square
\end{align*}

\emph{Proof of Lemma~\ref{lemma2}:~}
Similar to \eqref{not} we define 
\begin{eqnarray*}
\hat{T}^{(k)}_r (E_N,M,X) & = & \# \left\{ n_:~0 \leq n < M, ~n \equiv r ~\textup{mod}~ k, ~\textrm{and}~(e_{n+1},\dots,e_{n+k}) = X \right\} \nonumber
\end{eqnarray*}
and note that
$$
T (E_N,M,X) = \sum_{r=0}^{k-1} \hat{T}^{(k)}_r (E_N,M,X)
$$
Then, as in \eqref{fre}, we have
\begin{eqnarray}
& & \p \left( \max_{d < k \leq \log_2 N} ~\max_{X \in \{-1,1\}^k} ~\max_{1 \leq M \leq N+1-k} \left| T (E_N,M,X) - \frac{M}{2^k} \right| \geq \frac{16 \sqrt{N}}{d} \right) \nonumber\\
& \leq & \sum_{d < k \leq \log_2 N} ~\sum_{X \in \{-1,1\}^k} \p \left(\max_{1 \leq M \leq N+1-k}  \left| T(E_N,M,X) - \frac{M}{2^k} \right| \geq \frac{16 \sqrt{N}}{d} \right) \nonumber\\
& \leq & \sum_{d < k \leq \log_2 N} ~\sum_{X \in \{-1,1\}^k} \sum_{r=0}^{k-1} \p \left(\max_{1 \leq M \leq N+1-k}  \left| \hat{T}^{(k)}_r (E_N,M,X) - \frac{M}{k 2^k} \right| \geq \frac{16 \sqrt{N}}{dk} \right). \label{this}
\end{eqnarray}
Note that 
$$
\max_{1 \leq M \leq N+1-k}  \left| \hat{T}^{(k)}_r (E_N,M,X) - \frac{M}{k 2^k} \right| \leq \frac{1}{2^k} + \max_{\substack{1 \leq M \leq N+1-k,\\(M-r)/k \in \Z}}  \left| \hat{T}^{(k)}_r (E_N,M,X) - \frac{M}{k 2^k} \right|
$$
Again the functions $\hat{T}^{(k)}_r (E_N,M,X)$ are constructed in such a way that they are a sum of independent random variables, namely 
$$
\hat{T}^{(k)}_r (E_N,M,X) = \sum_{m=1}^{\lfloor (M-r)/k \rfloor} \mathds{1}_X \left( e_{(m-1)k+r+1}, \dots, e_{mk+r} \right).
$$
Furthermore, the random variables $\mathds{1}_X\left( e_{mk+r+1}, \dots, e_{m(k+1)+r} \right) - 2^{-k}$ have variance at most $2^{-k}$. Thus, using Lemma~\ref{bernstein}, we obtain for any $r$ and for $t = 3 \cdot 2^{-k/2} \sqrt{N} \sqrt{\log d}$, 
\begin{eqnarray*}
& & \p \left(\max_{\substack{1 \leq M \leq N+1-k,\\(M-r)/k \in \Z}}  \left| \hat{T}^{(k)}_r (E_N,M,X) - \frac{M}{k 2^k} \right| \geq 3 \cdot 2^{-k/2} \sqrt{N} \sqrt{\log d}\right) \\
& \leq & 2 \exp \left(- \frac{9 \cdot 2^{-k} N \log d}{2Nk^{-1}2^{-k} + 6 \cdot 2^{-k/2} \sqrt{N} \sqrt{\log d} ~d^{-1/2}/3} \right) \\
& \leq & 2 d^{-3k}
\end{eqnarray*}
for sufficiently large $N$. Note that for $k > d \geq 4$ we have
$$
\frac{1}{2^k} + 3 \cdot 2^{-k/2} \sqrt{\log d} \leq \frac{16}{dk},
$$
and consequently
$$
\p \left(\max_{1 \leq M \leq N+1-k}  \left| \hat{T}^{(k)}_r (E_N,M,X) - \frac{M}{k 2^k} \right| \geq \frac{16 \sqrt{N}}{dk} \right) \leq 2 d^{-3k}.
$$
Together with \eqref{this} this implies
\begin{eqnarray*}
& & \p \left( \max_{d < k \leq \log_2 N} ~\max_{X \in \{-1,1\}^k} ~\max_{1 \leq M \leq N+1-k} \left| T(E_N,M,X) - \frac{M}{2^k} \right| \geq \frac{16 \sqrt{N}}{d} \right) \\
& \leq & \sum_{d < k \leq \log_2 N} ~\sum_{X \in \{-1,1\}^k} \sum_{r=0}^{k-1} 2 d^{-3k} \\ 
& \leq & 2 \sum_{k=d+1}^\infty d^{-2k} \\
& \leq & \frac{1}{d^{2d}}. \qquad \square
\end{eqnarray*}

\section{Proof of Theorem~\ref{th1}}

To prove the existence of the limit distribution, we have to show the following: for any given $\ve>0$ and $t \in \R$ there exists a number $N_0=N_0(\ve)$ such that for $N_1,N_2 \geq N_0$ we have
\begin{equation} \label{cauchy}
\p \left( \no(E_{N_1}) \leq t \sqrt{N_1} \right) - \p \left( \no(E_{N_2}) \leq t \sqrt{N_2} \right) \leq \ve.
\end{equation}
By \eqref{n1} the limit distribution exists for $t=0$, and satisfies $F(0)=0$. Thus we can assume in the sequel that $t>0$. Let $\delta > 0$ be fixed (and ``small''), and choose $d$ sufficiently large such that all the inequalities $6 \sqrt{\log d}/\sqrt{d} \leq \delta/3, ~1/(d^2-1) \leq \ve/6,~16/d \leq \delta/3$ and $1/d^{2d} \leq \ve/6$ hold. By \eqref{ter1} and \eqref{ter2} we have
\begin{eqnarray}
& & \p \left(\no(E_{N_1}) \leq t \sqrt{N_1} \right) \label{arg1}\\
& \leq & \p \left(\nd(E_{N_1}) \leq (t+\delta) \sqrt{N_1} \right) 
\nonumber\\ 
& & + ~\underbrace{\p \left( \max_{k \leq d} ~\max_X ~\max_M \left| \bar{T}^{(d)}(E_{N_1},M,X) - \left(\frac{M+k-1}{d} -1 \right)\frac{k-1}{2^k}\right| > \delta/3 \sqrt{N_1} \right)}_{\leq \ve/6~ \textrm{by Lemma~\ref{lemma1}}} \nonumber\\
& & + ~\underbrace{\p \left( \max_{d < k \leq \log_2 N} ~\max_X ~\max_M \left| T(E_{N_1},M,X) - \frac{M}{2^k} \right| > \delta/3 \sqrt{N_1} \right)}_{\leq \ve/6~ \textrm{by Lemma~\ref{lemma2}}} \nonumber
\end{eqnarray}
for sufficiently large $N_1$. Similarly, by \eqref{ter3},
\begin{eqnarray*}
& & \p \left( \no(E_{N_2}) \leq t \sqrt{N_2} \right) \\
& \geq & \p \left( \nd(E_{N_2}) \leq (t-\delta) \sqrt{N_2} \right) \\
& & - ~\underbrace{\p \left( \max_{k \leq d} ~\max_X ~\max_M \left| \bar{T}^{(d)}(E_{N_2},M,X) - \left(\frac{M+k-1}{d} -1 \right)\frac{k-1}{2^k}\right| > \delta \sqrt{N_2} \right)}_{\leq \ve/6~ \textrm{by Lemma~\ref{lemma1}}},
\end{eqnarray*}
for sufficiently large $N_2$. By Lemma~\ref{lemmacru}
$$
\left| \p \left(\nd(E_{N_1}) \leq (t+\delta) \sqrt{N_1} \right) - \p \left( \frac{Z(s)}{\sqrt{d}} \in P^{(t+\delta)} ~\textup{for all~} s \in [0,1] \right) \right| \leq \ve/6
$$
and
$$
\left| \p \left(\nd(E_{N_2}) \leq (t-\delta) \sqrt{N_2} \right) - \p \left( \frac{Z(s)}{\sqrt{d}} \in P^{(t-\delta)} ~\textup{for all~} s \in [0,1] \right) \right| \leq \ve/6
$$
for sufficiently large $N_1$ and $N_2$, and since the exiting probabilities 
$$
\p \left( \frac{Z(s)}{\sqrt{d}} \in P^{(u)} ~\textup{for all~} s \in [0,1] \right)
$$
depend on the parameter $u$ continuously we have
$$
\p \left( \frac{Z(s)}{\sqrt{d}} \in P^{(t+\delta)} ~\textup{for all~} s \in [0,1] \right) - \p \left( \frac{Z(s)}{\sqrt{d}} \in P^{(t-\delta)} ~\textup{for all~} s \in [0,1] \right) \leq \ve/6
$$
if $\delta$ is sufficiently small. Overall, we have established \eqref{cauchy}, which proves the existence of the limit distribution of the normalized normality measure.\\

Now we turn to the continuity of the limit distribution function $F(t)$. Let $\ve>0$ be given. We have to show that
\begin{equation} \label{cont1}
F(t+\delta) - F(t) \leq \ve
\end{equation}
for some sufficiently small $\delta=\delta(\ve)$. The continuity of $F(t)$ at $t=0$ follows from \eqref{n1}. Thus we can henceforth assume that $t>0$. We have
\begin{eqnarray*}
&  & F(t+\delta) - F(t) \\
& = & \lim_{N \to \infty} \p \left( \no(E_N)/\sqrt{N} \leq t + \delta \right) - \lim_{N \to \infty} \p \left( \no(E_N)/\sqrt{N} \leq t \right).
\end{eqnarray*}
Arguing as in the lines following \eqref{arg1} we can show that
$$
\p \left( \no(E_N)/\sqrt{N} \leq t + \delta \right) \leq \p \left( \frac{Z(s)}{\sqrt{d}} \in P^{(t+2\delta)} ~\textup{for all~} s \in [0,1] \right) + \ve/3
$$
and
$$
\p \left( \no(E_N)/\sqrt{N} \leq t \right) \geq \p \left( \frac{Z(s)}{\sqrt{d}} \in P^{(t-\delta)} ~\textup{for all~} s \in [0,1] \right) - \ve/3
$$
for sufficiently large $N$. Together with the fact that
$$
\left| \p \left( \frac{Z(s)}{\sqrt{d}} \in P^{(t+2\delta)} ~\textup{for all~} s \in [0,1] \right) - \p \left( \frac{Z(s)}{\sqrt{d}} \in P^{(t-\delta)} ~\textup{for all~} s \in [0,1] \right) \right| \leq \ve/3
$$
for sufficiently small $\delta$, this establishes \eqref{cont1}. Altogether, we have proved Theorem~\ref{th1}.


\begin{thebibliography}{10}

\bibitem{ai}
C.~Aistleitner.
\newblock On the limit distribution of the well-distribution measure of random
  binary sequences.
\newblock {\em J. Theor. Nombres Bordeaux}, 25(2):245--259, 2013.

\bibitem{ai2}
C.~Aistleitner.
 \newblock Normal numbers and the normality measure. 
 \newblock \emph{Combin. Probab. Comput.}, 22(3):342--345, 2013. 

\bibitem{akmmr_min}
N.~Alon, Y.~Kohayakawa, C.~Mauduit, C.~G. Moreira, and V.~R{\"o}dl.
\newblock Measures of pseudorandomness for finite sequences: minimal values.
\newblock {\em Combin. Probab. Comput.}, 15(1-2):1--29, 2006.

\bibitem{akmmr_typ}
N.~Alon, Y.~Kohayakawa, C.~Mauduit, C.~G. Moreira, and V.~R{\"o}dl.
\newblock Measures of pseudorandomness for finite sequences: typical values.
\newblock {\em Proc. Lond. Math. Soc. (3)}, 95(3):778--812, 2007.

\bibitem{an}
V.~Anantharam.
\newblock A technique to study the correlation measures of binary sequences.
\newblock {\em Discrete Math.}, 308(24):6203--6209, 2008.

\bibitem{bpt}
I.~Berkes, W.~Philipp, and R.~F. Tichy.
\newblock Pseudorandom numbers and entropy conditions.
\newblock {\em J. Complexity}, 23(4-6):516--527, 2007.

\bibitem{cms}
J.~Cassaigne, C.~Mauduit, and A.~S{\'a}rk{\"o}zy.
\newblock On finite pseudorandom binary sequences. {VII}. {T}he measures of
  pseudorandomness.
\newblock {\em Acta Arith.}, 103(2):97--118, 2002.

\bibitem{einmahl}
U.~Einmahl and D.~M. Mason.
\newblock Some universal results on the behavior of increments of partial sums.
\newblock {\em Ann. Probab.}, 24(3):1388--1407, 1996.

\bibitem{g1}
K.~Gyarmati.
\newblock On the correlation of binary sequences.
\newblock {\em Studia Sci. Math. Hungar.}, 42(1):79--93, 2005.

\bibitem{g2}
K.~Gyarmati and C.~Mauduit.
\newblock On the correlation of binary sequences, {II}.
\newblock {\em Discrete Math.}, 312(5):811--818, 2012.

\bibitem{mat}
J.~Matou{\v{s}}ek.
\newblock {\em Geometric discrepancy}, volume~18 of {\em Algorithms and
  Combinatorics}.
\newblock Springer-Verlag, Berlin, 2010.
\newblock An illustrated guide, Revised paperback reprint of the 1999 original.

\bibitem{ms1}
C.~Mauduit and A.~S{\'a}rk{\"o}zy.
\newblock On finite pseudorandom binary sequences. {I}. {M}easure of
  pseudorandomness, the {L}egendre symbol.
\newblock {\em Acta Arith.}, 82(4):365--377, 1997.

\bibitem{mer}
L.~M{\'e}rai.
\newblock Construction of pseudorandom binary sequences over elliptic curves
  using multiplicative characters.
\newblock {\em Publ. Math. Debrecen}, 80(1-2):199--213, 2012.

\bibitem{rivat}
J.~Rivat and A.~S{\'a}rk{\"o}zy.
\newblock Modular constructions of pseudorandom binary sequences with composite
  moduli.
\newblock {\em Period. Math. Hungar.}, 51(2):75--107, 2005.

\bibitem{schmidt}
K.~Schmidt.
\newblock The correlation measures of finite sequences: limiting distributions and minimum values.
\newblock Preprint.  Available at \url{http://arxiv.org/abs/1404.0172}.

\bibitem{sz}
B.~Sziklai.
\newblock On the symmetry of finite pseudorandom binary sequences.
\newblock {\em Unif. Distrib. Theory}, 6(2):143--156, 2011.

\bibitem{ww}
W.~Whitt.
\newblock {\em Stochastic-process limits}.
\newblock Springer Series in Operations Research. Springer-Verlag, New York,
  2002.

\bibitem{wu}
C.~Wu, X.~Weng, and Z.~Chen.
\newblock Construction of {$k$}-ary pseudorandom elliptic curve sequences.
\newblock {\em Wuhan Univ. J. Nat. Sci.}, 16(5):452--456, 2011.

\end{thebibliography}
\end{document}